\documentclass [12pt] {article}

\usepackage{latexsym}

\font\bb=msbm10 at 12pt
\newcommand{\R}{\mbox{{\bb R}}}
\newcommand{\D}{\mbox{{\bb D}}}
\newcommand{\C}{\mbox{{\bb C}}}

\everymath{\displaystyle}

\newenvironment{multi}
{ \begin{equation} \begin{array}{c} }
{ \end{array} \end{equation} }

\makeatletter

\@addtoreset{equation}{section}
\makeatother

\begin{document}

\title{Lie  derivatives and structure Jacobi operator on real hypersurfaces in complex projective spaces II}
\author{Juan de Dios P\'{e}rez and David P\'{e}rez-L\'{o}pez}
\date{}
\maketitle

\begin{abstract} Let $M$ be a real hypersurface in complex projective space. The almost contact metric structure on $M$ allows us to consider, for any nonnull real number $k$, the corresponding $k$-th generalized Tanaka-Webster connection on $M$ and, associated to it, a differential operator of first order of Lie type. Considering such a differential operator and Lie derivative we define, from the structure Jacobi operator $R_{\xi}$ on $M$ a tensor field of type (1,2), $R_{{\xi}_T}^{(k)}$. We obtain some classifications of real hypersurfaces for which $R_{{\xi}_T}^{(k)}$ is either symmetric or skew symmetric.
\end{abstract}



\section{Introduction.}
Consider a complex projective space $\C P^m$, $m\geq 2$, endowed with the complex structure $J$ and the Fubini-Study metric $g$ of constant holomorphic sectional curvature $4$. Let $M$ be a {\it connected real hypersurface} of $\C P^m$ without boundary, $\nabla$ its Levi-Civita connection on $M$ and $N$ a local unit normal vector field on $M$. Then $\xi =-JN$ is a tangent vector field to $M$ called the structure vector field on $M$ (or the Reeb vector field). The Kaehlerian structure $(J,g)$ of $\C P^m$ induces on $M$ an almost contact metric structure $(\phi , \xi , \eta , g)$, see \cite{B},  where $\phi$ is the tangent component of $J$, $\eta$ is an 1-form given by $\eta(X)=g(X,\xi)$ for any $X$ tangent to $M$ and $g$ is the metric induced on $M$.                                                                                                               

Takagi, see  \cite{T1}, \cite{T2}, \cite{T3}, classified homogeneous real hypersurfaces in $\C P^m$ into 6 types. Kimura, \cite{K1}, proved that these 6 types are the unique Hopf real hypersurfaces in $\C P^m$ with constant principal curvatures. A real hypersurface $M$ is called Hopf if the structure vector field is principal, that is, $A\xi =\alpha\xi$ for a certain function $\alpha$ on $M$, where $A$ denotes the shape operator of $M$ associated to $N$. We will also denote by $\D$ the maximal holomorphic distribution on $M$, given by all vector fields orthogonal to $\xi$.                          

 Among the real hypersurfaces appearing in Takagi's list we mention the following ones:

Type $(A_1)$ real hypersurfaces are geodesic hyperspheres of radius $r$, $0 < r < \frac {\pi}{2}$. They have 2 distinct constant principal curvatures, $2\cot 2r$ with eigenspace $\R [\xi]$ and $\cot r$ with eigenspace $\D$.

Type $(A_2)$ are tubes of radius $r$, $0 < r < \frac {\pi}{2}$, over totally geodesic complex projective spaces $\C P^n$, $0 < n < m-1$. They have 3 distinct constant principal curvatures, $2\cot 2r$ with eigenspace $\R [\xi]$, $\cot r$ and $-\tan r$. The corresponding eigenspaces of $\cot r$ and $-\tan r$ are complementary and $\phi$-invariant distributions in $\D$.

From now on we will call type $(A)$ real hypersurfaces to both of either type $(A_1)$ or type $(A_2)$. The other real hypersurfaces appearing in Takagi's list have either 3 or 5 distinct constant principal curvatures.

Tanaka, \cite{TA}, and, independtly Webster, \cite{W}, studied the canonical affine connection defined on a non-degenerate, pseudo-Hermitian CR-manifold, known as the Tanaka-Webster connection. Tanno, \cite{TAN},  generalized this connection for contact metric manifolds. 

Cho, using the almost contact metric structure on a real hypersurface $M$ of $\C P^m$, generalized Tanno's connection, defining, for any non-zero real number $k$,  the $k$-th generalized Tanaka-Webster connection $\hat{\nabla}^{(k)}$, see \cite{CH1}, \cite{CH2}, by

\begin{multi} \label{1.1}
\hat{\nabla}^{(k)}_XY=\nabla_XY+g(\phi AX,Y)\xi -\eta (Y)\phi AX-k\eta (X)\phi Y
\end{multi}

\noindent for any X,Y tangent to $M$. Then $\hat{\nabla}^{(k)}\eta =0$, $\hat{\nabla}^{(k)}\xi =0$, $\hat{\nabla}^{(k)}g=0$, $\hat{\nabla}^{(k)}\phi =0$. In particular, if the shape operator of a real hypersurface satisfies $\phi A+A \phi =2k\phi$, $M$ is a contact manifold and the $k$-th generalized Tanaka-Webster connection coincides with the Tanaka-Webster connection. From (\ref{1.1}) the $k$-th Cho operator on $M$ associated to a tangent vector field $X$ is defined by $F_X^{(k)}Y=g(\phi AX,Y)\xi -\eta(Y)\phi AX-k\eta(X)\phi Y$, for any $Y$ tangent to $M$. The torsion of the $k$-th generalized Tanaka-Webster connection is then given by $T^{(k)}(X,Y)=F^{(k)}_XY-F^{(k)}_YX$. For any $X$ tangent to $M$ we call $T_X^{(k)}$ to the operator on $M$ given by $T_X^{(k)}Y=T^{(k)}(X,Y)$, for any $Y$ tangent to $M$ and call it the torsion operator associated to $X$.

Jacobi fields along geodesics of a given Riemannian manifold $(\tilde{M},\tilde{g})$ satisfy a very well-known differential equation. This classical differential equation naturally inspires the so-called {\it Jacobi operator}. That is, if $\tilde{R}$ is the curvature operator of $\tilde{M}$, and $X$ is any tangent vector field to $\tilde{M}$, the Jacobi operator (with respect to $X$) at $p\in M$, $\tilde{R}_X\in$End$(T_p\tilde{M})$, is defined as $(\tilde{R}_XY)(p)=(\tilde{R}(Y,X)X)(p)$ for all $Y\in T_p\tilde{M}$, being a selfadjoint endomorphism of the tangent bundle $T\tilde{M}$ of $\tilde{M}$. Clearly, each tangent vector field $X$ to $\tilde{M}$ provides a Jacobi operator with respect to $X$. 

If $M$ is a real hypersurface of $\C P^{m}$ the Jacobi operator associated to the structure vector field $\xi$ is called the structure Jacobi operator on $M$ and we denote it by $R_{\xi}$.

Let  $\cal L$ denote the Lie derivative on $M$. Therefore ${\cal L}_XY=\nabla_XY-\nabla_YX$ for any $X,Y$ tangent to $M$. If we consider $\hat{\nabla}^{(k)}$, we can define a differential operator of first order ${\cal L}^{(k)}$ on $M$ by ${\cal L}_X^{(k)}Y=\hat{\nabla}_X^{(k)}Y-\hat{\nabla}_Y^{(k)}X={\cal L}_XY+T_X^{(k)}Y$, for any $X,Y$ tangent to $M$.

We can associate to $R_{\xi}$ a tensor field of type (1,2), $R_{{\xi}_T}^{(k)}$, given by $R_{{\xi}_T}^{(k)}(X,Y)=\lbrack T_X^{(k)},R_{\xi} \rbrack Y=T_X^{(k)}R_{\xi}Y-R_{\xi}T_X^{(k)}Y$, for any $X,Y$ tangent to $M$.

In the first paper of this series, \cite{P}, first author studied the condition ${\cal L}^{(k)}R_{\xi}={\cal L}R_{\xi}$, that is, $({\cal L}_X^{(k)}R_{\xi})Y=({\cal L}_XR_{\xi})Y$ for any $X,Y$ tangent to $M$. This condition is equivalent to the fact that $R_{{\xi}_T}^{(k)}$ vanishes, obtaining the following

\noindent {\bf Theorem} {\it There do not exist real hypersurfaces in $\C P^m$, $m \geq 3$, such that ${\cal L}^{(k)}R_{\xi}={\cal L}R_{\xi}$.}

The purpose of this paper is to study weaker conditions on $M$. Indeed, we will consider real hypersurfaces $M$ in $\C P^m$ whose structure Jacobi operator satisfies the following Codazzi condition

\begin{multi} \label{1.2}
(({\cal L}_X^{(k)}-{\cal L}_X)R_{\xi})Y=(({\cal L}_Y^{(k)}-{\cal L}_Y)R_{\xi})X
\end{multi}

\noindent for any $X,Y$ tangent to $M$. It is easy to see that this condition is equivalent to $R_{{\xi}_T}^{(k)}$ being symmetric, that is, $R_{{\xi}_T}^{(k)}(X,Y)=R_{{\xi}_T}^{(k)}(Y,X)$ for any $X,Y$ tangent to $M$. We will prove the following theorems.

\noindent {\bf Theorem 1} {\it There does not exist any Hopf real hypersurface $M$ in $\C P^m$, $m \geq 3$, such that $R_{\xi}$ satisfies (\ref{1.2}) for some nonnull real number $k$.}

\noindent and

\noindent {\bf Theorem 2} {\it There does not exist any non-Hopf real hypersurface $M$ in $\C P^m$, $m \geq 3$, such that $R_{\xi}$ satisfies (\ref{1.2}) for any nonnull real number $k$ such that $kg(A\xi ,\xi) \neq -1$.}

A different weaker condition is a Killing type condition given by

\begin{multi} \label{1.3}
(({\cal L}_X^{(k)}-{\cal L}_X)R_{\xi})Y+(({\cal L}_Y^{(k)}-{\cal L}_Y)R_{\xi})X=0
\end{multi}

\noindent for any $X,Y$ tangent to $M$. Now this condition is equivalent to the fact that $R_{{\xi}_T}^{(k)}$ is skew symmetric, that is, $R_{{\xi}_T}^{(k)}(X,Y)=-R_{{\xi}_T}^{(k)}(Y,X)$ for any $X,Y$ tangent to $M$. In this case we have a quite different situation as given by the following results.

\noindent {\bf Theorem 3} {\it Let $M$ be a Hopf real hypersurface in $\C P^m$, $m \geq 3$, and $k$ a nonnull real number. Then $M$ satisfies (\ref{1.3}) if and only if either $k=1$ and $M$ is locally congruent to a geodesic hypersphere of radius $\frac{\pi}{4}$ or $k^2 \geq 4$ and $M$ is locally congruent to a tube of radius $r$, $0 < r <\frac{\pi}{2}$, around $\C P^n$, $0 < n < m-1$, whose radius satisfies either $cot(r)=k$ or $cot(r)=-\frac{1}{k}$.} 

In the non-Hopf case we obtain

\noindent {\bf Theorem 4} {\it There does not exist any non-Hopf real hypersurface $M$ in $\C P^m$, $m \geq 3$, such that $R_{{\xi}_T}^{(k)}$ is skew symmetric for any nonnull real number $k$ satisfying $kg(A\xi ,\xi)=1$.}

\noindent {\bf Aknowledgements} This work was supported by MINECO-FEDER Project  MTM 2016-78807-C2-1-P.

\section{Preliminaries.}

Throughout this paper, all manifolds, vector fields, etc., will be considered of class $C^{\infty}$ unless otherwise stated. Let $M$ be a connected real hypersurface in $\C P^m$, $m\geq 2$, without boundary. Let $N$ be a locally defined unit normal vector field on $M$. Let $\nabla$ be the Levi-Civita connection on $M$ and $(J,g)$ the Kaehlerian structure of $\C P^m$.

For any vector field $X$ tangent to $M$ we write $JX=\phi X+\eta(X)N$, where $\phi X$ denotes the tangential component of $JX$ and $-JN=\xi$. Then $(\phi,\xi,\eta,g)$ is an almost contact metric structure on $M$, see \cite{B}. That is, we have

\begin{equation} \label{2.1}
\phi^2X=-X+\eta(X)\xi,  \quad  \eta(\xi)=1,  \quad   g(\phi X,\phi Y)=g(X,Y)-\eta(X)\eta(Y)
\end{equation}

\noindent for any tangent vectors $X,Y$ to $M$. From (\ref{2.1}) we obtain

\begin{multi} \label{elem}
\phi\xi =0,   \quad   \eta(X)=g(X,\xi).
\end{multi}

\noindent From the parallelism of $J$ we get

\begin{equation} \label{lasphi}
(\nabla_X\phi)Y=\eta(Y)AX-g(AX,Y)\xi
\end{equation}

\noindent and

\begin{equation} \label{elem2}
\nabla_X\xi =\phi AX
\end{equation}

\noindent for any $X,Y$ tangent to $M$, where $A$ denotes the shape operator of the immersion. As the ambient space has holomorphic sectional curvature $4$, the equations of Gauss and Codazzi are given, respectively, by

\begin{equation} \label{Gauss}
\begin{array}{rl}
R(X,Y)Z =  \displaystyle g(Y,Z)X - g(X,Z)Y + g(\phi Y,Z)\phi X - g(\phi X,Z)\phi Y      \\
- 2g(\phi X,Y)\phi Z + g(AY,Z)AX - g(AX,Z)AY     \\
\end{array}\end{equation}

\noindent and

\begin{equation} \label{Codazzi}
(\nabla_XA)Y-(\nabla_YA)X=\eta(X)\phi Y-\eta(Y)\phi X-2g(\phi X,Y)\xi
\end{equation}

\noindent for any vector fields $X,Y, Z$ tangent to $M$, where $R$ is the curvature tensor of $M$. We will call the maximal holomorphic distribution $\D$ on $M$ to the following one: at any $p \in M$, $\D (p)=\{ X\in T_pM \vert  g(X,\xi)=0\}$. If $U \in \D $ we define ${\D}_U=\{ X \in \D \vert g(X,U)=g(X, \phi U)=0 \}$. We will say that $M$ is Hopf if $\xi$ is principal, that is, $A\xi =\alpha\xi$ for a certain function $\alpha$ on $M$.

From Gauss equation the structure Jacobi operator on $M$ is given by

\begin{equation}\label{perez}
R_{\xi}X=X-\eta(X)\xi +g(A\xi ,\xi)AX-g(A\xi ,X)A\xi
\end{equation}

\noindent for any $X$ tangent to $M$

In the sequel we need the following results:

\noindent {\bf Theorem 2.1. \cite{M}} {\it
Let $M$ be a Hopf real hypersurface of $\C P^m$, $m \geq 2$, and let $X \in \D$ such that $AX=\lambda X$. Then $\alpha =g(A\xi ,\xi)$ is constant, $2\lambda -\alpha \neq 0$ and $\phi X$ is principal with principal curvature $\frac {\alpha\lambda +2}{2\lambda -\alpha}$.}

\noindent {\bf Theorem 2.2. \cite{O}} {\it
Let $M$ be a real hypersurface of $\C P^m$, $m \geq 2$. Then the following are equivalent:
\begin{enumerate}
\item $M$ is locally congruent to a real hypersurface of type $(A)$.
\item $\phi A=A\phi$.
\end{enumerate}  }

\noindent {\bf Theorem 2.3 \cite{PS}} {\it
There do not exist real hypersurfaces $M$ in $\C P^m$, $m \geq 3$, whose shape operator is given by $A\xi =\alpha\xi +\beta U$, $AU=\beta\xi +\frac{\beta^2-1}{\alpha} U$, $A\phi U=-\frac{1}{\alpha} \phi U$, where $U$ is a unit vector field in $\D$, $\alpha$ and $\beta$ are nonvanishing functions on $M$, the eigenvalues of $A$ in ${\D}_U$ are different from $0$, $-\frac{1}{\alpha}$ and $\frac{\beta^2-1}{\alpha}$ and if $Z \in {\D}_U$ satisfies $AZ=\lambda Z$, then $A\phi Z=\lambda \phi Z$.}

\noindent {\bf Theorem 2.4. \cite{PS}} {\it
There do not exist real hypersurfaces $M$ in ${\C} P^m$, $m \geq 3$, whose shape operator is given by $A\xi =\xi +\beta U$, $AU=\beta\xi +(\beta^2-1)U$, $A\phi U=-\phi U$, for a unit $U \in {\D}$, being $\beta$ a nonvanishing function, and there exists $Z \in {\D}_U$ such that $AZ=-Z$, $A\phi Z=-\phi Z$.}

\section{Proofs of Theorems 1 and 2.}

Let $M$ be a real hypersurface in $\C P^m$ satisfying (\ref{1.2}). Then $F_X^{(k)}R_{\xi}Y-F_{R_{\xi}Y}^{(k)}X-R_{\xi}F_X^{(k)}Y+R_{\xi}F_Y^{(k)}X=F_Y^{(k)}R_{\xi}X-F_{R_{\xi}X}^{(k)}Y-R_{\xi}F_Y^{(k)}X+R_{\xi}F_X^{(k)}Y$ for any $X,Y$ tangent to $M$. This yields                                            

\begin{multi} \label{3.1}
g(\phi AX,R_{\xi}Y)\xi -k\eta(X)\phi R_{\xi}Y-g(\phi AR_{\xi}Y,X)\xi +\eta(X)\phi AR_{\xi}Y+\eta(Y)R_{\xi}\phi AX  \\
+k\eta(X)R_{\xi}\phi Y-\eta(X)R_{\xi}\phi AY-k\eta(Y)R_{\xi}\phi X=g(\phi AY,R_{\xi}X)\xi -k\eta(Y)\phi R_{\xi}X \\
-g(\phi AR_{\xi}X,Y)\xi +\eta(Y)\phi AR_{\xi}X+\eta(X)R_{\xi}\phi AY+k\eta(Y)R_{\xi}\phi X-\eta(Y)R_{\xi}\phi AX     \\  -k\eta(X)R_{\xi}\phi Y  \\
\end{multi}

\noindent for any $X,Y$ tangent to $M$.

If in (\ref{3.1}) we take $X,Y \in \D$ we obtain

\begin{multi} \label{3.2}
g(\phi AX,R_{\xi}Y)\xi-g(\phi AR_{\xi}Y,X)\xi =g(\phi AY,R_{\xi}X)\xi -g(\phi AR_{\xi}X,Y)\xi
\end{multi}

\noindent for any $X,Y \in \D$. And if we take $X=\xi$, $Y \in \D$ in (\ref{3.1}) it follows

\begin{multi} \label{3.3}
g(\phi A\xi ,R_{\xi}Y)\xi -k\phi R_{\xi}Y+\phi AR_{\xi}Y+2kR_{\xi}\phi Y-2R_{\xi}\phi AY=0
\end{multi}

\noindent for any $Y \in \D$.

Let us suppose that $M$ is Hopf with $A\xi =\alpha\xi$. Then (\ref{3.3}) gives

\begin{multi} \label{3.4}
-k\phi R_{\xi}Y+\phi AR_{\xi}Y+2kR_{\xi}\phi Y-2R_{\xi}\phi AY=0
\end{multi}

\noindent for any $Y \in \D$. If $X \in \D$ then $R_{\xi}X=X+\alpha AX$. Therefore, (\ref{3.2}) yields $g(\phi AX,Y+\alpha AY)-g(\phi AY+\alpha\phi A^2Y,X)=g(\phi AY,X+\alpha AX)-g(\phi AX+\alpha\phi A^2X,Y$ for any $X,Y \in \D$. Then $2g(\phi AX,Y)+2\alpha g(\phi AX,AY)-2g(\phi AY,X)-\alpha g(\phi A^2Y,X)+\alpha g(\phi A^2X,Y)=0$, for any $X,Y \in \D$. As $M$ is Hopf this implies

\begin{multi} \label{3.5}
2\phi AX+2\alpha A\phi AX+2A\phi X+\alpha A^2\phi X+\alpha \phi A^2X=0
\end{multi}

\noindent for any $X \in \D$. Let now $X \in \D$ such that $AX=\lambda X$. Then $A\phi X=\mu\phi X$, $\mu =\frac{\alpha\lambda +2}{2\lambda -\alpha}$. Then (\ref{3.5}) yields $2\lambda +2\alpha\lambda\mu +2\mu +\alpha\mu^2+\alpha\lambda^2=0$. That is, $2(\lambda +\mu)+\alpha(\lambda +\mu)^2=0$. Thus if we suppose $\lambda +\mu =0$, $\frac{2\lambda^2+2}{2\lambda -\alpha}=0$, which is impossible.  Therefore,

\begin{multi} \label{3.6}
\alpha \neq 0
\end{multi}

\noindent and $\alpha(\lambda +\mu)=-2$. That is, $\frac{2\lambda^2+2}{2\lambda -\alpha}=-\frac{2}{\alpha}$. This implies $(\lambda^2+1)\alpha =\alpha -2\lambda$, that is, $\lambda^2\alpha =-2\lambda$. Therefore, either $\lambda =0$ or $\lambda =-\frac{2}{\alpha}$. If $X$ satisfies $AX=0$, then $A\phi X=-\frac{2}{\alpha}\phi X$. This implies $R_{\xi}X=X$ and $R_{\xi}\phi X=\phi X+\alpha A\phi X=-\phi X$. From (\ref{3.4}) we get $-k\phi X-2k\phi X=0$, that is, $3k\phi X=0$, which is impossible. We obtain a similar situation if we suppose that $AX=-\frac{2}{\alpha} X$ and Theorem 1 is proved.

Suppose now that $M$ is non-Hopf. Then, at least on a neighbourhood of a certain point of $M$ we can find a unit vector field $U \in \D$ and functions $\alpha$ and $\beta$, where $\beta$ does not vanish, such that we can write $A\xi =\alpha\xi +\beta U$. From now on we consider all the calculations on such a neighborhood. Then (\ref{3.3}) gives

\begin{multi} \label{3.7}                                                                                                                                        \beta g(\phi U,R_{\xi}Y)\xi -k\phi R_{\xi}Y +\phi AR_{\xi}Y+2kR_{\xi}\phi Y-2R_{\xi}\phi AY=0
\end{multi}

\noindent for any $Y \in \D$. Its scalar product with $\xi$ yields $\beta g(R_{\xi}\phi U,Y)=0$ for any $Y \in \D$. As $R_{\xi}\xi =0$ we obtain $R_{\xi}\phi U=0$ and then

\begin{multi} \label{3.8}
A\phi U=-\frac{1}{\alpha}\phi U.
\end{multi}

If in (\ref{3.7}) we take $Y=\phi U$ we get $-2kR_{\xi}U+\frac{2}{\alpha} R_{\xi}\phi^2U=0$. Thus $-(2k+\frac{2}{\alpha})R_{\xi}U=0$ and, as we suppose that $\alpha \neq -\frac{1}{k}$ we obtain $R_{\xi}U=0$ or

\begin{multi} \label{3.9}
AU=\beta\xi +\frac{\beta^2-1}{\alpha} U.
\end{multi}

Therefore, $\D_U$ is $A$-invariant. From (\ref{3.2}) we know that $g((\phi A+A\phi)X,R_{\xi}Y)=g((\phi A+A\phi)Y,R_{\xi}X)$ for any $X,Y \in \D_U$. Therefore,

\begin{multi} \label{3.10}
-(\phi A+A\phi)R_{\xi}Y=R_{\xi}(\phi A+A\phi)Y
\end{multi}

\noindent for any $Y \in \D_U$. If we take $Y \in \D_U$ unit and such that $AY=\lambda Y$ in (\ref{3.7}) we have $-k\phi R_{\xi}Y+\phi AR_{\xi}Y+2(k-\lambda)R_{\xi}\phi Y=0$. Then, $-(1+\alpha\lambda)k\phi Y+(1+\alpha\lambda)\phi AY+2(k-\lambda)(\phi Y+\alpha A\phi Y)=0$. That is, $((1+\alpha\lambda)(\lambda -k)+2(k-\lambda))\phi Y+2\alpha(k-\lambda)A\phi Y=0$. Thus $2\alpha(k-\lambda)A\phi Y=(\alpha\lambda -1)(k-\lambda)\phi Y$. Therefore, either $k=\lambda$ or if $\lambda \neq k$,

\begin{multi} \label{3.11}
A\phi Y=\frac{\alpha\lambda -1}{2\alpha}\phi Y.
\end{multi}

We will call $\mu =\frac{\alpha\lambda -1}{2\alpha}$. If in (\ref{3.7}) we take $\phi Y$ instead of $Y$ we have $-k\phi R_{\xi}\phi Y+\phi AR_{\xi}\phi Y-2kR_{\xi}Y-2R_{\xi}\phi A\phi Y=0$. This implies $-k(1+\alpha\mu)\phi^2Y+(1+\alpha\mu)\phi A\phi Y-2k(1+\alpha\lambda)Y+2\mu (1+\alpha\lambda)Y=0$. Then $(-1+\alpha\mu -2\alpha\lambda)(k-\mu)=0$. This gives, bearing in mind the value of $\mu$, $-3(1+\alpha\lambda)(k-\mu)=0$. Therefore, we have two possibilities: either $\lambda =-\frac{1}{\alpha}$ and then $\mu =-\frac{1}{\alpha}$, or $\mu =k$ and, as we suppose $\lambda \neq k$, in this case $\alpha\lambda -1=2k\alpha$ and $\lambda =\frac{2\alpha k+1}{\alpha}$.

As $-(\phi A+A\phi)R_{\xi}Y=R_{\xi}(\phi A+A\phi)Y$, we have $-(1+\alpha\lambda)(\phi A+A\phi)Y=(\lambda +\mu)R_{\xi}\phi Y$. Then $-(1+\alpha\lambda)(\lambda +\mu)=(1+\alpha\mu)(\lambda +\mu)$. If $\lambda =\mu =-\frac{1}{\alpha}$, we obtain $\lambda +\mu \neq 0$.

If $\mu =k$  and we suppose $\lambda +\mu \neq 0$, it follows $-1-\alpha\lambda =1+\alpha\mu$ and $\lambda +\mu =-\frac{2}{\alpha}$. Then $\frac{2\alpha k+1}{\alpha} +k=\frac{3\alpha k+1}{\alpha}=-\frac{2}{\alpha}$ yields $3\alpha k=-3$ and this contradicts the fact of $\alpha k \neq -1$. Therefore, $\lambda +\mu =0$ and $3\alpha k+1=0$ yields $\alpha =-\frac{1}{3k}$. In this case, $\alpha$, $\lambda$ and $\mu$ are constant.

If the unique principal curvature on $\D_U$ is $k$, as $-(1+\alpha k)(\phi A+A\phi)Y=2k(1+\alpha k)\phi Y$, we should have $4k(1+\alpha k)=0$ and, as $k \neq 0$, $\alpha k=-1$, which is impossible. Therefore, on $\D_U$ there must be, at least, a principal curvature $\lambda$ different from $k$.

Let us suppose that on $\D_U$ there is a principal curvature $\lambda =-\frac{1}{\alpha}$ and let $Y$ be a unit eigenvector corresponding to $\lambda$. As we have seen $A\phi Y=-\frac{1}{\alpha}\phi Y$ and the equation of Codazzi implies $(\nabla_YA)\phi Y-(\nabla_{\phi Y}A)Y=-2\xi$. This gives $-Y(\frac{1}{\alpha})\phi Y-\frac{1}{\alpha}\nabla_Y\phi Y-A\nabla_Y\phi Y+(\phi Y)(\frac{1}{\alpha})Y+\frac{1}{\alpha}\nabla_{\phi Y}Y+A\nabla_{\phi Y}Y=-2\xi$. If we take its scalar product with $\xi$ we obtain $\frac{1}{\alpha} g(\phi Y,\phi AY)- g(\nabla_Y\phi Y,\alpha\xi +\beta U)-\frac{1}{\alpha} g(Y,\phi A\phi Y)+g(\nabla_{\phi Y}Y,\alpha\xi +\beta U)=-2$. This implies

\begin{multi} \label{3.12}
\beta g(\lbrack \phi Y,Y \rbrack ,U)=\frac{2}{\alpha^2}.
\end{multi}

And its scalar product with $U$ yields $-\frac{1}{\alpha} g(\nabla_Y\phi Y,U)+\frac{1}{\alpha} g(\nabla_{\phi Y}Y,U)-g(\nabla_Y\phi Y,\beta\xi +\frac{\beta^2-1}{\alpha} U)+g(\nabla_{\phi Y}Y,\beta\xi +\frac{\beta^2-1}{\alpha} U)=0$. This gives $\frac{\beta^2}{\alpha} g(\lbrack \phi Y,Y \rbrack ,U)-\frac{2\beta}{\alpha} =0$. Thus

\begin{multi} \label{3.13}
\beta g(\lbrack \phi Y,Y \rbrack ,U)=2.
\end{multi}

Fron (\ref{3.12}) and (\ref{3.13}) we have $\alpha^2 =1$ and, taking $-\xi$ instead of $\xi$ if necessary, we can suppose $\alpha =1$. From Theorem 2.4 this kind of real hypersurfaces does not exist.

Now the unique principal curvatures on $\D_U$ are $k$ and $-k$, with $A\xi =-\frac{1}{3k}\xi +\beta U$, $AU=\beta\xi +3k(1-\beta^2)U$, $A\phi U=3k\phi U$. Let a unit $X \in \D_U$ such that $AX=\lambda X$, where we know that $\lambda$ is constant. Then we have $(\nabla_XA)\xi -(\nabla_{\xi}A)X=-\phi X$, that is, $\nabla_X(-\frac{1}{3k}\xi +\beta U)-A\phi AX-\nabla_{\xi}(\lambda X)+A\nabla_{\xi}X=-\phi X$. Its scalar product with $\xi$ implies $g(\nabla_{\xi}X,-\frac{1}{3k}\xi +\beta U)=0$. As we suppose $\beta \neq 0$, we obtain $g(\nabla_{\xi}X,U)=0$. But its scalar product with $U$ gives $X(\beta)-\lambda g(\nabla_{\xi}X,U)+g(\nabla_{\xi}X,\beta\xi +3k(\beta^2-1)U)=0$. Then, $X(\beta)+(3k(1-\beta^2)-\lambda)g(\nabla_{\xi}X,U)=0$. this yields

\begin{multi} \label{3.14}
X(\beta)=0
\end{multi}

\noindent for any $X \in \D_U$.

On the other hand $(\nabla_UA)\xi -(\nabla_{\xi}A)U=-\phi U$ gives $-\frac{1}{3k}\phi AU+U(\beta)U+\beta\nabla_UU-A\phi AU-\xi(\beta)\xi -\beta\phi A\xi-3k\xi(1-\beta^2)U-3k(1-\beta^2)\nabla_{\xi}U+A\nabla_{\xi}U=\phi U$. Its scalar product with $\xi$ implies $-\xi(\beta)+g(\nabla_{\xi}U,-\frac{1}{3k}\xi +\beta U)=0$. Therefore,

\begin{multi} \label{3.15}
\xi(\beta)=0.
\end{multi}

From the scalar product of the above expression and $U$ we get $U(\beta)-3k\xi(1-\beta^2)+g(\nabla_{\xi}U,\beta\xi +3k(1-\beta^2)U)=0$. From (\ref{3.15}) we have

\begin{multi} \label{3.16}
U(\beta)=0.
\end{multi}

Once more, Codazzi equation implies $(\nabla_{\phi U}A)\xi -(\nabla_{\xi}A)\phi U=U$, that is, $-\frac{1}{3k}\phi A\phi U+(\phi U)(\beta)U+\beta\nabla_{\phi U}U+3kAU-3k\nabla_{\xi}\phi U+A\nabla_{\xi}\phi U=U$. Its scalar product with $\xi$ gives $-\beta g(U,\phi A\phi U)+3k\beta +3kg(\phi U,\phi A\xi)+g(\nabla_{\xi}\phi U,-\frac{1}{3k}\xi +\beta U)=0$. Thus $9k\beta +\frac{1}{3k} g(\phi U,\phi A\xi)+\beta g(\nabla_{\xi}\phi U,U)=0$. As $\beta \neq 0$ we obtain

\begin{multi} \label{3.17}
g(\nabla_{\xi}U,\phi U)=9k+\frac{1}{3k}.
\end{multi}

Its scalar product with $U$ yields $\frac{1}{3k} g(A\phi U,\phi U)+(\phi U)(\beta)+9k^2(1-\beta^2)-3kg(\nabla_{\xi}\phi U,U)+g(\nabla_{\xi}\phi U,\beta\xi +3k(1-\beta^2)U)=1$. From (\ref{3.17}) we obtain $(\phi U)(\beta)+9k^2-9k^2\beta^2+27k^2+1-\beta^2+3k(\beta^2-1)(9k+\frac{1}{3k})=0$ and this implies

\begin{multi} \label{3.18}
(\phi U)(\beta)=-9k^2(1+2\beta^2).
\end{multi}

From (\ref{3.14}), (\ref{3.15}), (\ref{3.16}) and (\ref{3.18}) we obtain

\begin{multi} \label{3.19}
grad(\beta)=\omega\phi U
\end{multi}

\noindent where $\omega =-9k^2(1+2\beta^2)$.

As $g(\nabla_Xgrad(\beta),Y)=g(\nabla_Ygrad(\beta),X)$, for any $X,Y$ tangent to $M$, we obtain $X(\omega)g(\phi U,Y)+\omega g(\nabla_X\phi U,Y)=Y(\omega)g(\phi U,X)+\omega g(\nabla_Y\phi U,X)$, for any $X,Y$ tangent to $M$. Taking $Y=\xi$ we get $\omega g(\nabla_X\phi U,\xi)=\omega g(\nabla_{\xi}\phi U,X)$, that is, $-\omega g(U,AX)=\omega g(\nabla_{\xi}\phi U,X)$ for any $X$ tangent to $M$. If $X=U$ we have $-\omega g(AU,U)=\omega g(\nabla_{\xi}\phi U,U)=-\omega g(\nabla_{\xi}U,\phi U)$. Therefore, either $\omega =0$, which is impossible, or $g(AU,U)=g(\nabla_{\xi}U,\phi U)$. This yields $3k(1-\beta^2)=9k+\frac{1}{3k}$, that is, $1-\beta^2=3+\frac{1}{9k^2}$. Equivalently, $\beta^2=-2-\frac{1}{9k^2}$, which is impossible. This finishes the proof of Theorem 2.

\section{Proofs of Theorems 3 and 4.}

If $M$ is a real hypersurface satisfying (\ref{1.3}) we obtain

\begin{multi} \label{4.1}
g(\phi AX,R_{\xi}Y)\xi -k\eta(X)\phi R_{\xi}Y-g(\phi AR_{\xi}Y,X)\xi +\eta(X)\phi AR_{\xi}Y   \\
+\eta(Y)R_{\xi}\phi AX+k\eta(X)R_{\xi}\phi Y-\eta(X)R_{\xi}\phi AY-k\eta(Y)R_{\xi}\phi X+g(\phi AY,R_{\xi}X)\xi    \\
-k\eta(Y)\phi R_{\xi}X-g(\phi AR_{\xi}X,Y)\xi +\eta(Y)\phi AR_{\xi}X+\eta(X)R_{\xi}\phi AY   \\
+k\eta(Y)R_{\xi}\phi X-\eta(Y)R_{\xi}\phi AX-k\eta(X)R_{\xi}\phi Y=0      \\
\end{multi}

\noindent for any $X,Y$ tangent to $M$. If in (\ref{4.1}) we take $X,Y \in \D$ we have

\begin{multi} \label{4.2}
g(\phi AX,R_{\xi}Y)\xi -g(\phi AR_{\xi}Y,X)\xi +g(\phi AY,R_{\xi}X)\xi -g(\phi AR_{\xi}X,Y)\xi =0
\end{multi}

\noindent for any $X,Y \in \D$. Taking $X=\xi$, $Y \in \D$ in (\ref{4.1}) we get $g(\phi A\xi ,R_{\xi}Y)\xi -k\phi R_{\xi}Y+\phi AR_{\xi}Y+kR_{\xi}\phi Y-R_{\xi}\phi AY+R_{\xi}\phi AY-kR_{\xi}\phi Y=0$. This gives

\begin{multi} \label{4.3}
g(\phi A\xi ,R_{\xi}Y)\xi -k\phi R_{\xi}Y+\phi AR_{\xi}Y=0
\end{multi}

\noindent for any $Y \in \D$. 

Let us suppose that $M$ is Hopf and $A\xi =\alpha\xi$. From (\ref{4.2}) we get $-A\phi R_{\xi}Y-\phi AR_{\xi}Y+R_{\xi}\phi AY+R_{\xi}A\phi Y=0$ for any $Y \in \D$. Take a unit $Y \in \D$ such that $AY=\lambda Y$. Then $A\phi Y=\mu\phi Y$, $\mu =\frac{\alpha\lambda +2}{2\lambda -\alpha}$. In this case, $R_{\xi}Y=(1+\alpha\lambda)Y$ and $R_{\xi}\phi Y=(1+\alpha\mu)\phi Y$. From (\ref{4.3}) we have $-k(1+\alpha\lambda)\phi Y+\lambda(1+\alpha\lambda)\phi Y=0$. Thus $(\lambda -k)(1+\alpha\lambda)=0$. Therefore, if $\alpha =0$, $\lambda =k$ and if $\alpha \neq 0$, either $\lambda =k$ or $\lambda =-\frac{1}{\alpha}$. But as $-A\phi R_{\xi}Y-\phi AR_{\xi}Y+R_{\xi}\phi AY+R_{\xi}A\phi Y=0$ we get $-(1+\alpha\lambda)\mu\phi Y-(1+\alpha\lambda)\lambda\phi Y-\lambda(1+\alpha\mu)\phi Y+\mu(1+\alpha\mu)\phi Y=0$. That is, $(\lambda +\mu)(\mu -\lambda)\alpha =0$. As $\lambda +\mu =\lambda +\frac{\alpha\lambda +2}{2\lambda -\alpha} =\frac{2\lambda^2+2}{2\lambda -\alpha} \neq 0$, it follows $\alpha(\mu -\lambda)=0$.

If $\alpha =0$ the unique principal curvature on $\D$ is $k$. From \cite{CR} $M$ must be locally congruent to a geodesic hypersphere of radius $\frac{\pi}{4}$. Therefore $k=cot(\frac{\pi}{4})=1$.

If $\alpha \neq 0$, $\lambda =\mu$, that is $\phi A=A\phi$ and from Theorem 2.2, $M$ must be locally congruent to a real hypersurface of type $(A)$. Its principal curvatures are $\alpha$, $-\frac{1}{\alpha}$ and $k$. In the case of a geodesic hypersphere we should have $\alpha =-\frac{1}{k}$. Thus $2cot(2r)=-\frac{1}{k}$ and $cot(r)=k$. Thus $-\frac{1}{k} =2cot(2r)=cot(r)-tan(r)=k-\frac{1}{k}$ yields $k=0$, which is impossible.

Therefore, $M$ is locally congruent to a tube of radius $r$, $0 < r < \frac{\pi}{4}$, around $\C P^n$, $0 < n < m-1$. Their principal curvatures are $\alpha =2cot(2r)=cot(r)-tan(r)$, $cot(r)$ and $-tan(r)$. We thus have two situations: either  $cot(r)=k$ and $tan(r)=\frac{1}{\alpha}$, $\alpha =k-\frac{1}{\alpha}$, or $cot(r)=-\frac{1}{\alpha}$, $tan(r)=-k$ and again $\alpha =-\frac{1}{\alpha} +k$. Anyway, $\alpha^2-k\alpha +1=0$, Thus $\alpha =\frac{k \pm \sqrt{k^2-4}}{2}$. Therefore, $k^2 \geq 4$ and either $cot(r)=k$ or $cot(r)=-\frac{1}{k}$, finishing the proof of Theorem 3.

Let us suppose now that $M$ is non-Hopf and write, as in previous section, $A\xi =\alpha\xi +\beta U$. Then (\ref{4.3}) gives

\begin{multi} \label{4.4}
\beta g(\phi U,R_{\xi}Y)\xi -k\phi R_{\xi}Y+\phi AR_{\xi}Y=0
\end{multi}

\noindent for any $Y \in \D$. The scalar product of (\ref{4.4}) and $\xi$ yields $\beta g(R_{\xi}\phi U,Y)=0$ for any $Y \in \D$ and, being $\beta \neq 0$ and $R_{\xi}\xi =0$, we obtain $R_{\xi}\phi U=0$, that is, $\phi U+\alpha A\phi U=0$. Then

\begin{multi} \label{4.5}
\alpha \neq 0
\end{multi}

\noindent and

\begin{multi} \label{4.6}
A\phi U=-\frac{1}{\alpha}\phi U. 
\end{multi}

Also (\ref{4.4}) becomes

\begin{multi} \label{4.7}
\phi AR_{\xi}Y=k\phi R_{\xi}Y
\end{multi}

\noindent for any $Y \in \D$. If we take (\ref{4.7}) into (\ref{4.2}) we obtain $g(\phi AX,R_{\xi}Y)-kg(\phi R_{\xi}Y,X)+g(\phi AY,R_{\xi}X)-kg(\phi R_{\xi}X,Y)=0$ for any $X,Y \in \D$. As $g(R_{\xi}\phi AX+kR_{\xi}\phi X-A\phi R_{\xi}X-k\phi R_{\xi}X,\xi)=0$ because $R_{\xi}\phi U=0$, we get

\begin{multi} \label{4.8}
R_{\xi}\phi AX+kR_{\xi}\phi X-A\phi R_{\xi}X-k\phi R_{\xi}X=0
\end{multi}

\noindent for any $X \in \D$. Taking $X=\phi U$ in (\ref{4.8}) we have $R_{\xi}\phi A\phi U-kR_{\xi}U=0$. From (\ref{4.6}) it follows $(\frac{1}{\alpha} -k)R_{\xi}U=0$. As we suppose $k\alpha \neq 1$ we obtain $R_{\xi}U=0$ and then

\begin{multi} \label{4.9}
AU=\beta\xi +\frac{\beta^2-1}{\alpha} U.
\end{multi}

Now $\D_U$ is A-invariant. Take a unit $Y \in \D_U$ such that $AY=\lambda Y$. Then $R_{\xi}Y=(1+\alpha\lambda)Y$ and applying $\phi$ to (\ref{4.7}) we have $-AR_{\xi}Y+g(AR_{\xi}Y,\xi)\xi =-kR_{\xi}Y+kg(R_{\xi}Y,\xi)\xi$. As $R_{\xi}U=0$ this yields $AR_{\xi}Y=kR_{\xi}Y$. Therefore $\lambda(1+\alpha\lambda)Y=k(1+\alpha\lambda)Y$. Thus either $\lambda =-\frac{1}{\alpha}$ or $\lambda =k$.

Moreover, from (\ref{4.8}), for such a $Y$ we get $\lambda R_{\xi}\phi Y+kR_{\xi}\phi Y-(1+\alpha\lambda)A\phi Y-(1+\alpha\lambda)k\phi Y=0$. Now if $\lambda =-\frac{1}{\alpha}$ we arrive to $(k-\frac{1}{\alpha})R_{\xi}\phi Y=0$ and, being $k\alpha \neq 1$, we obtain $R_{\xi}\phi Y=0=\phi Y+\alpha A\phi Y$, which implies $A\phi Y=-\frac{1}{\alpha}\phi Y$.

If $\lambda =k$, from (\ref{4.8}) we have $2kR_{\xi}\phi Y-(1+k\alpha)A\phi Y-k(1+k\alpha)\phi Y=0$, that is, $2k(\phi Y+\alpha A\phi Y)-(1+k\alpha)A\phi Y-k(1+k\alpha)\phi Y=0$. Thus $(k\alpha -1)A\phi Y=k(k\alpha -1)\phi Y$ and, as $k\alpha \neq 1$ we obtain $A\phi Y=k\phi Y$. We have obtained that any eigenspace in $\D_U$ is $\phi$-invariant.

From Theorem 2.4 the case of $\lambda =-\frac{1}{\alpha}$ does not occur and then the unique principal curvature on $\D_U$ is $k$. Then, for any $X \in \D_U$, Codazzi equation gives $(\nabla_XA)\phi X-(\nabla_{\phi X}A)X=-2\xi$. This yields $k\nabla_X\phi X-A\nabla_X\phi X-k\nabla_{\phi X}X+A\nabla_{\phi X}X=-2\xi$. Its scalar product with $\xi$ implies $-kg(\phi X,\phi AX)-g(\nabla_X\phi X,\alpha\xi +\beta U)+kg(X,\phi A\phi X)+g(\nabla_{\phi X}X,\alpha\xi +\beta U)=-2$, that is

\begin{multi} \label{4.10}
\beta g(\lbrack \phi X,X \rbrack ,U)=2k^2-2k\alpha -2.
\end{multi}

Its scalar product with $U$ gives $-kg(\lbrack \phi X,X \rbrack ,U)-g(\nabla_X\phi X,\beta\xi +\frac{\beta^2-1}{\alpha} U)+g(\nabla_{\phi X}X,\beta\xi +\frac{\beta^2-1}{\alpha} U)=0$. This yields

\begin{multi} \label{4.11}
(\frac{\beta^2-1}{\alpha} -k)g(\lbrack \phi X,X \rbrack ,U)=-2k\beta.
\end{multi}

If $k=\frac{\beta^2-1}{\alpha}$, $2\beta k=0$, which is impossible. Then $k \neq 0$ and $k\neq \frac{\beta^2-1}{\alpha}$. Moreover, from (\ref{4.10}) and (\ref{4.11}) we obtain $2(k^2-k\alpha -1)(\frac{\beta^2-1}{\alpha}-k) =-2k\beta^2$. That is, $(k^2-(k\alpha +1))(\beta^2-(k\alpha +1))=-k\alpha\beta^2$. If $k=-\frac{1}{\alpha}$, we should have $k^2\beta^2=\beta^2$ and therefore, $k^2=1$ and $\alpha^2=1$. This hypersurfaces do not exist from Theorem 2.4.

Thus the unique principal curvature on $\D_U$ is $k \neq 0,-\frac{1}{\alpha},\frac{\beta^2-1}{\alpha}$ and such a real hypersurface does not exist from Theore, 2.3. This finishes our proof.

 \hfill $\Box$

{\sc
\begin{trivlist}

\item Juan de Dios Perez: jdperez@ugr.es  \\
Departamento de Geometria y Topologia  and IEMATH\\
Universidad de Granada \\
18071 Granada  \\
Spain  \\

\item David Perez-Lopez: davidpl109@correo.ugr.es  \\
Fundaci\'{o}n I+D del Software Libre - FIDESOL\\
Avda. de la Innovaci\'{o}n, 1 - Ed. BIC CEEI (PTS)\\
18016 Armilla (Granada)\\
Spain\\

\end{trivlist}
}

\begin{thebibliography}{99}

\bibitem{B} D.E. Blair, Riemannian Geometry of contact and symplectic manifolds, {\em Progress in Mathematics} {\bf 203} (2002), Birkhauser Boston Inc. Boston.

\bibitem{CH1} J.T. Cho, CR-structures on real hypersurfaces of a complex space form, {\em Publ. Math. Debrecen} {\bf 54} (1999), 473-487.

\bibitem{CR} T.E. Cecil and P.J. Ryan, Focal sets and real hypersurfaces in complex projective space, {\em Trans. A.M.S.} {\bf 269} (1982), 481-499.

\bibitem{CH2} J.T. Cho, Pseudo-Einstein CR-structures on real hypersurfaces in a complex space form, {\em Hokkaido Math. J.} {\bf 37} (2008), 1-17.

\bibitem{JPS} I. Jeong, E. Pak and Y.J. Suh, {\em Real hypersurfaces in complex two-plane Grassmannians with generalized Tanaka-Webster invariant shape operator}, {\em J. Math. Phys. Anal. Geom.} {\bf 9} (2013), 360-378.

\bibitem{K1} M. Kimura, Real hypersurfaces and complex submanifolds in complex projective space, {\em Trans. A.M.S.} {\bf 296} (1986), 137-149.





\bibitem{M} Y. Maeda, On real hypersurfaces of a complex projective space, {\em J. Math. Soc. Japan} {\bf 28} (1976), 529.540.

\bibitem{O} M. Okumura, On some real hypersurfaces of a complex projective space, {\em Trans. A.M.S.} {\bf 212} (1975), 355-364.

\bibitem{P} J.D. P\'{e}rez, Lie derivatives and structure Jacobi operator on real hypersurfaces in complex projective spaces, {\em Diff. Geom. Appl.} {\bf 50} (2017), 1-10.

\bibitem{PS} J.D. P\'{e}rez and F.G. Santos, Real hypersurfaces in complex projective space whose structure Jacobi operator satisfies ${\cal L}_{\xi}R_{\xi}=\nabla_{\xi}R_{\xi}$, {\em Rocky Mount. J. Math.} {\bf 39} (2009), 1293-1301.



\bibitem{T1} R. Takagi, On homogeneous real hypersurfaces in a complex projective space, {\em Osaka J. Math.} {\bf 10} (1973), 495-506.

\bibitem{T2} R. Takagi, Real hypersurfaces in complex projective space with constant principal curvatures, {\em J. Math. Soc. Japan} {\bf 27} (1975), 43-53.

\bibitem{T3} R. Takagi, Real hypersurfaces in complex projective space with constant principal curvatures II, {\em J. Math. Soc. Japan} {\bf 27} (1975), 507-516.

\bibitem{TA} N. Tanaka, On non-degenerate real hypersurfaces, graded Lie algebras and Cartan connections, {\em Japan. J. Math.} {\bf 2} (1976), 131-190.

\bibitem{TAN} S. Tanno, Variational problems on contact Riemennian manifolds, {\em Trans. A.M.S.} {\bf 314} (1989), 349-379.

\bibitem{W} S.M. Webster, Pseudohermitian structures on a real hypersurface, {\em J. Diff. Geom.} {\bf 13} (1978), 25-41.

\end{thebibliography}
\end{document}